\documentclass[a4paper,11pt,oneside]{amsart}
\usepackage{mathptmx}       
\usepackage{helvet}         
\usepackage{courier}        
\usepackage{type1cm}        
\usepackage{graphicx}        
\usepackage{multicol}        
\usepackage[bottom]{footmisc}
\usepackage{url}

\usepackage{mathtools}
\DeclarePairedDelimiter{\floor}{\lfloor}{\rfloor}

\usepackage{amsfonts}
\usepackage{amsmath}
\usepackage{amssymb}

\begin{document}

\title{The inverse problem for the fractional conductivity equation: a survey}

\author{Giovanni Covi}
\address{Department of Mathematics and Statistics, University of Helsinki, Pietari Kalmin katu 5 (Exactum), Helsinki, Finland}
\email{giovanni.covi@helsinki.fi}

\begin{abstract}
    The fractional Calderón problem asks to determine the unknown coefficients in a nonlocal, elliptic equation of fractional order from exterior measurements of its solutions. There has been substantial work on many aspects of this inverse problem. In this review we collect some recent results related to the conductivity formulation of the fractional Calderón problem.
\end{abstract} 
\maketitle
\vspace{15mm}

This short note discusses the results achieved in recent years in the study of the fractional Calderón problem in its conductivity formulation. We start by briefly discussing the fractional Calderón problem in the form in which it was introduced in \cite{GSU16}, see also the survey \cite{Sal17} for more information on the topic. Let $s\in (0,1)$, $n\in \mathbb N$, and consider the fractional Laplace operator $(-\Delta)^s$ defined via the Fourier transform $\mathcal F$ as 
$$(-\Delta)^su = \mathcal F^{-1}(|\xi|^{2s}\mathcal F u(\xi)),$$
or equivalently via a singular integral as 
$$ (-\Delta)^su(x) := c_{n,s}PV\int_{\mathbb R^n}\frac{u(x)-u(y)}{|x-y|^{n+2s}}dy, $$
where $c_{n,s}:=\frac{4^s\Gamma(n/2+s)}{\pi^{n/2}|\Gamma(-s)|}$ is a fixed constant. Let $\Omega$ be a bounded open set, and consider the following fractional Schrödinger equation
$$\begin{cases}
    (-\Delta)^su+qu =0, & \Omega \\ u =f, & \Omega_e:=\mathbb R^n\setminus\overline\Omega
\end{cases},$$
where $q$ is a fixed unknown potential and $f$ is an exterior value belonging to the Sobolev space $H^s(\Omega_e)$. In this formulation, the fractional Calderón problem consists in recovering $q$ from the nonlocal Dirichlet-to-Neumann (DN) map $\Lambda_q$ defined as
$$ \Lambda_q : H^s(\Omega_e)\ni f \mapsto (-\Delta)^su|_{\Omega_e}\in H^s(\Omega_e)^*. $$

This inverse problem was solved in \cite{GSU16} for $q\in L^\infty$, and later in \cite{RS17b} for $q\in L^{n/2s}$. Uniqueness and reconstruction were achieved even for a single measurement in \cite{GRSU18}. Exponential stability results in line with the classical case were proved in \cite{RS17a, RS17b, Rül20}, while Lipschitz stability was proved in the case of finite dimension \cite{RS18}. The instability properties of the fractional Calderón problem were explored in \cite{BCR24} using the technique developed in \cite{KRS21}. Many variations of the fractional Schrödinger equation have been studied, such as with general local (\cite{CMRU22}, see also \cite{CLR18, CMR21}), quasilocal (\cite{Cov21}) and magnetic perturbations (\cite{Cov20b, Li19, Li20a, Li20b}), or on closed (\cite{FKU24}) and complete Riemannian manifolds (\cite{CO23}). Some related nonlocal operators were considered in e.g. \cite{BS22, CdHS24, CGR22, GLX17, KLW21, LLR19}. Finally, the relations between the classical and fractional Calderón problems were object of study in \cite{CGRU23, GU21, Rül23}. We refer to the survey \cite{Sal17} for more (somewhat older) references on the fractional Calderón problem in the Schrödinger formulation. The conductivity formulation of the fractional Calderón problem was introduced in \cite{Cov20a}, building on the nonlocal vector calculus developed in \cite{DGLZ12, DGLZ13}. Let $s\in (0,1)$, $n\in \mathbb N$, and consider the fractional gradient $\nabla^s$ defined as
$$ \nabla^su(x,y):= \frac{c^{1/2}_{n,s}}{\sqrt{2}}\frac{u(y)-u(x)}{|x-y|^{n/2+s+1}}(x-y) $$
for a smooth and compactly supported scalar function $u$. Since $\nabla^s$ is a bounded operator (\cite[eq. (7)]{Cov20a}), it is possible to extend it by density to act as $\nabla^s: H^s(\mathbb R^n)\rightarrow L^2(\mathbb R^{2n})$, and to observe (\cite[Lemma 4.1]{Cov20a}) that it holds $\lim_{s\rightarrow 1^-}\|\nabla^su\|_{L^2(\mathbb R^{2n})} = \|\nabla u\|_{L^2(\mathbb R^{n})}$.
The fractional divergence $(\nabla\cdot)^s:L^2(\mathbb R^{2n})\rightarrow H^{-s}(\mathbb R^{n})$ is defined as the adjoint of the fractional gradient, i.e. by the relation $$ \langle (\nabla\cdot)^sv, u \rangle_{H^{-s}(\mathbb R^{n}), H^s(\mathbb R^n)} := \langle v, \nabla^su \rangle_{L^2(\mathbb R^{2n})},  $$
and it verifies the fundamental identity $(\nabla\cdot)^s\nabla^s = (-\Delta)^s$ involving the fractional Laplacian (\cite[Lemma 2.1]{Cov20a}). If now $\gamma\in L^\infty(\mathbb R^n)$ is a scalar conductivity, it makes sense to consider the operator $\mathbf C^s_\gamma : H^s(\mathbb R^n)\rightarrow H^{-s}(\mathbb R^n)$ defined by
$$ \mathbf C^s_\gamma u := (\nabla\cdot)^s(\Theta\cdot\nabla^s u), $$
where $\Theta(x,y):= \gamma^{1/2}(x)\gamma^{1/2}(y)I$. This is the (isotropic) fractional conductivity operator. The inverse problem for the fractional conductivity equation 
$$\begin{cases}
    \mathbf C^s_\gamma u = (\nabla\cdot)^s(\Theta\cdot\nabla^s u) =0, & \Omega \\ u =f, & \Omega_e
\end{cases},$$
consists in recovering the conductivity $\gamma$ from the nonlocal DN map 
$$ \Lambda_\gamma : H^s(\Omega_e)\ni f \mapsto \mathbf C^s_\gamma u|_{\Omega_e}\in H^s(\Omega_e)^*. $$
In the rest of the review we will briefly discuss the results achieved for this inverse problem and the techniques involved in their proofs.

\section{Fractional random walks}
The fractional conductivity operator arises naturally from the study of random walks with long jumps on graphs via a homogenization procedure. Let $h>0$, and consider a random walk on the lattice $h\mathbb Z^n$ with time steps $h^{2s}\mathbb N$ and transition probability
$$ P(x,y)=\frac{1}{m(x)}\frac{\gamma^{1/2}(x)\gamma^{1/2}(y)}{|x-y|^{n+2s}}, \qquad x\neq y, \quad x,y\in h\mathbb Z^n, $$
where $m(x)$ is a normalization factor. Let $u(x,t)$ be the probability of finding the random walker at $x\in h\mathbb Z^n$ at time $t\in h^{2s}\mathbb N$. Then one can prove (\cite[Sec. 5]{Cov20a}, see also \cite{Val09} for the case $\gamma\equiv 1$) that in the limit $h\rightarrow 0$ the function $u$ solves 
$$ \gamma^{1/2}(x)\partial_t u(x,t) \approx \mathbf C^s_\gamma u(x,t), $$
so that the equation $\mathbf C^s_\gamma u=0$ corresponds to the equilibrium state $\partial_tu\equiv 0$. The recent work \cite{CL24} generalizes the random walk described above to a general graph $G$, and studies the inverse problem of recovering both the graph $G$ and the conductivity $\gamma$ from partial random walk data. Using a new algebraic method, \cite[Th. 1.2]{CL24} shows that the recovery is possible up to natural constraints, thus opening a new line of research related to fractional conductivity on Riemannian manifolds \cite[Sec. 8]{CL24}. 

\section{The fractional Liouville reduction}
A familiar method for the study of the classical conductivity operator $\nabla\cdot\gamma\nabla$ is the \emph{Liouville reduction}, a change of variables transforming it into the Schrödinger operator $\Delta+q_\gamma$. A similar reduction holds in the fractional case (\cite[Th. 3.1]{Cov20a}): we have
\begin{equation*}
   \begin{cases}
    \mathbf C^s_\gamma u =0, & \Omega \\ u =f, & \Omega_e
\end{cases} \qquad\Leftrightarrow\qquad \begin{cases}
    ((-\Delta)^s+q_\gamma) w =0, & \Omega \\ w =\gamma^{1/2}f, & \Omega_e
\end{cases}, 
\end{equation*}
where $w=\gamma^{1/2}u$ and $q_\gamma = -\gamma^{-1/2}(-\Delta)^s\gamma^{1/2}$. Since the main part of $(-\Delta)^s+q_\gamma$ is independent of $\gamma$, and many powerful techniques (such as Runge approximation) are available for the transformed problem, this reduction greatly simplifies the study of fractional conductivity. As a result, uniqueness and reconstruction (even with a single measurement) were obtained in \cite[Th. 1.1, 1.2]{Cov20a}.

\section{Generalizations of the uniqueness results}

\noindent\textbf{Low regularity.} The regularity assumptions of \cite[Th 1.1]{Cov20a} are not sharp. The work \cite{RZ22a} considers the fractional conductivity equation in the case of low regularity, obtaining uniqueness for conductivities belonging to the Sobolev space $H^{s,n/s}$. Moreover, the proof is obtained by different means, without making use of the Runge approximation property of the fractional Laplacian. Quite interestingly, some counterexamples to uniqueness arise in the low regularity case.\\

\noindent\textbf{Magnetic perturbations.} The uniqueness properties of the inverse problem related to fractional conductivity with magnetic perturbations were studied in the work \cite[Th. 1.1]{Cov20b}, where it was shown that, similarly to the classical case, a natural gauge involving the magnetic and electric potentials $A,q$ exists. The studied operator is
$$   (\nabla^s + A)^*(\nabla^s + A) + q, \qquad\qquad  A\in L^2(\mathbb R^{2n}),\quad q\in L^\infty(\mathbb R^n), $$
which is related to a random walk model with long jumps and weights.\\

\noindent\textbf{The higher order problem.} The definition of the fractional gradient operator $\nabla^s$ was extended in \cite{CMR21} to include all the cases $s\in \mathbb R^+\setminus \mathbb Z$ in the following natural way:
$$ \nabla^s := \nabla^{s-\floor{s}}\nabla^{\floor{s}}, \qquad\qquad \floor{s}:=\max\{z\in\mathbb Z, z\leq s\}. $$
Using this, the authors were able to consider the higher order fractional conductivity equation (even with magnetic perturbations), and obtain uniqueness for the related inverse problem. This required also the high order fractional Poincar\'e inequality $$ \|(-\Delta)^{t/2}u\|_{L^2(\mathbb R^n)} \leq c_{n,s,K} \|(-\Delta)^{s/2}u\|_{L^2(\mathbb R^n)}, \qquad t\in[0,s], \; K \mbox{ compact}, \; u\in H^s_K(\mathbb R^n)  ,$$ which was obtained in the same work using a variety of equivalent proofs. \\

\noindent\textbf{The anisotropic case.} There is of course great interest in the study of the anisotropic version of the fractional Calder\'on problem in conductivity formulation. This was first considered in \cite{Cov22}, where a uniqueness result for the operator $$ \mathbf C^s_A:=(\nabla\cdot)^s(A\cdot\nabla^s) $$ was obtained for a large class of anisotropic conductivities $A$. The cited work uses a special matrix decomposition into separable matrix functions and a new reduction technique which was introduced in \cite{CdHS24} for the study of fractional elasticity.  \\

\noindent\textbf{Unbounded domains and the global problem.} All the results cited above make the (somewhat) reasonable assumption that the conductivity $\gamma$ is constant in $\Omega_e:=\mathbb R^n\setminus\overline\Omega$. The work \cite{RZ22b} is the first to drop this assumption, thus considering the recovery of a general scalar conductivity $\gamma$ in sets bounded only in one direction. The companion work \cite{RZ22c} by the same authors shows interesting counterexamples to uniqueness in this new framework. The so called \emph{global} inverse problem of fractional conductivity continues along this line of research, considering the recovery of the conductivity $\gamma$ on all of $\mathbb R^n$. This result is reminiscent of the classical boundary determination method by Kohn and Vogelius, in that a special class of concentrating solutions is constructed and used for obtaining exterior determination.

\section{Stability}
The stability properties of the inverse problem for the fractional conductivity equation were studied in \cite{CRTZ22}. The main result of this work is a logarithmic stability estimate in bounded domains, which holds under suitable a priori assumptions on the conductivity. These observations are further complemented by an instability result in the spirit of Mandache, which uses the technique of \cite{KRS21} in order to obtain the optimality of the logarithmic estimate. This is in line with the results known in both the classical and the fractional Schrödinger cases. It should however be noted that the methods used differ substantially at many technical points.





\end{document}